\newcommand{\hp}[1]{\hspace*{#1ex}}

\newcommand{\hh}{\hp1}

\newcommand{\parbf}[1]{\par\textbf{#1}}
\newcommand{\pbfni}[1]{\par\noindent\textbf{#1}}

\newcommand{\spar}{\\\hspace*{1em}}

\newcommand{\nrm}{\normalfont}

\newcommand{\cur}{^\mathsf{C}}
\newcommand{\unc}{^\text{\reflectbox{$\mathsf C$}}}
\newcommand{\cart}{\times}
\newcommand{\tCart}{\Pi}
\newcommand{\tDisu}{\Sigma}
\newcommand{\Cart}{\tCart\hp{0.2}}
\newcommand{\Disu}{\tDisu\hp{0.2}}
\newcommand{\insec}{\cap}
\newcommand{\union}{\cup}
\newcommand{\tUnion}{\raisebox{-.1 ex}{\large\pacmat{\cup}}}
\newcommand{\Union}{\tUnion\,}
\newcommand{\tInsec}{\raisebox{-.1 ex}{\large\pacmat{\cap}}}
\newcommand{\Insec}{\tInsec\,}
\newcommand{\tpll}{{\pacmat{{|}{|}}}}
\newcommand{\pll}{\,\tpll\,}
\newcommand{\Pll}{\text{\large${|}{|}$}}
\newcommand{\bool}{\mathbb{B}}

\newcommand{\isrel}{\mathrm{isrel}}
\newcommand{\tdom}{\mathcal{D}}
\newcommand{\dom}{\tdom\,}

\newcommand{\tran}{\mathcal{R}}
\newcommand{\ran}{\mathcal{R}\,}

\newcommand{\rna}[1]{\{#1\}}
\newcommand{\tcod}{\mathcal{C}}

\newcommand{\codt}{\mathrm{cod}\,}
\newcommand{\tgrph}{\mathcal{G}}

\newcommand{\real}{\mathbb{R}}

\newcommand{\ab}{\,.\,}
\newcommand{\ba}{\,|\,}

\newcommand{\arr}{\rightarrow}
\newcommand{\onto}{\rightarrowtriangle}
\newcommand{\evry}{\forall\,}

\newcommand{\membr}{\in}
\newcommand{\cmps}{\circ}
\newcommand{\bin}{:}
\newcommand{\sbset}{\subseteq}
\newcommand{\trsp}{^{\mathsf{T}}}
\newcommand{\mpty}{\emptyset}
\newcommand{\sor}{\vee}
\newcommand{\sand}{\wedge}
\newcommand{\siff}{\equiv}

\newcommand{\simp}{\Rightarrow}

\newlength{\wa}
\newlength{\wb} 
\newcommand{\overstrike} [2] {\text {%
\settowidth{\wa}{#1}\settowidth {\wb} {#2}%
#1\hskip-0.5\wa\hskip-0.5\wb #2\hskip0.5\wa\hskip-0.5\wb}}
\newcommand{\pacmat} [1] {\text{${#1}$}}
\newcommand{\tparr}{\overstrike{\pacmat{\rightarrow}}{\raisebox{0.25ex}{\pacmat{\scriptscriptstyle/\;}}}}

\newcommand{\parr}{\,\tparr\,}
\newcommand{\tdarr}{\raisebox{0.5ex}{\small...}\hspace{-.6ex}\raisebox{0.21ex}{$\scriptstyle\succ$}}

\newcommand{\tsrc}{\mathrm{src}}
\newcommand{\src}{\tsrc\,}
\newcommand{\ttgt}{\mathrm{tgt}}
\newcommand{\tgt}{\ttgt\,}
\newcommand{\pr}{\uppi}

\newcommand{\ide}{\mathrm{id}}
\newcommand{\tlrr}{\pacmat{\leftrightarrow\hspace{-1.8ex}\rightarrow}}
\newcommand{\tlngr}[1]{\stackrel{#1}{\longrightarrow}}

\newcommand{\lrar}{\tlrr}

\newcommand{\cnv}{^{\scriptscriptstyle\smile}}
\newcommand{\gam}{\gamma}

\newcommand{\udel}{\updelta}
\newcommand{\ulam}{\uplambda}
\documentclass{elsarticle}
\usepackage{amsmath,amssymb,stmaryrd,hyperref,upgreek,graphics}
\input{diagxy}

\newtheorem{defin}{Definition}
\newcommand{\dref}{Definition~\ref}
\newcommand{\drefs}{Definitions~\ref}
\newtheorem{therm}{Theorem}
\newcommand{\tref}{Theorem~\ref}

\begin{document}
\begin{frontmatter}
\title{Why mathematics needs engineering}
\author{Raymond Boute, INTEC, Ghent University \hp2 \texttt{raymond.boute@pandora.be}} 
%
%
\begin{abstract}
\noindent  Engineering needs mathematics, but the converse is also increasingly evident. Indeed,  mathematics is still recovering from the drawbacks of several ``reforms''.  Encouraging is the revived interest in proofs indicated by various recent {\em introduction to proof}-type textbooks.  Yet, many of these texts defeat their own purpose by self-conflicting definitions. Most affected are fundamental concepts such as relations and functions, despite flawless accounts 50 years ago. We take the viewpoint that definitions and theorems are tools for capturing, analyzing and understanding mathematical concepts and hence, like any tools, require diligent engineering.  This is illustrated for relations and functions, their algebraic properties and their relation to category theory, with the {\em Halmos principle} for definitions and the {\em Arnold principle} for axiomatization as design guidelines.
\end{abstract}
\end{frontmatter}
%
%
\section{Introduction: Mathematics and Engineering}\label{intro}
%
%
\noindent Mathematics has been intertwined with engineering since antiquity ~\cite{Boye91hima,Russ04fore}.
   \par Kline notes that ``{\em More than anything else mathematics is a method}\/''~\cite{Klin53miwc}.  Arguably, the primary purpose of this method is {\em effective reasoning}. This view best explains what Wigner calls {\em the unreasonable effectiveness of mathematics}~\cite{Wign60ueom}, in particular its practical usefulness far beyond the originally intended application areas.  From this perspective, the dichotomy between Platonism and formalism dissolves:  mathematical objects {\em do}\/ exist,  albeit in an abstract universe.  Formalism, definitions and theorems are the tools to study them.
   \par Tools, being artifacts, deserve careful design, borrowing criteria and guidelines from engineering.  Some of these also been discussed by Jos\'e Oliveira \cite{Oliv97form} in another context.  Here we focus on using engineering principles in mathematics.
   \par Foremost is enhancing the effectiveness in reasoning. Symbolic notation properly designed and used yields extra guidance via the shape of the expressions.  It should function like well-meshed gears in a Swiss precision clockwork.
   \par Aptness and economy in capturing the abstract objects of interest ensures conceptual malleability, generality and practical usefulness. Human factors are influential here, and it is often overlooked that this is a highly individual matter of temperament and background.  Even so, everyone benefits from clear conceptualization and reasoning.  For instance, {\em separation of concerns} avoids the common misconceptions caused by intellectual noise and conceptual tangling.  
   \par In classical mathematics, methods and notations were often thought-out carefully.  In algebra, for instance, symbolic notation started with Diophantus and evolved into its current form via Vi\`ete and Descartes~\cite{Bash99lcvd,Boye91hima}, rarely violating good design practices, thus making symbolic calculation today's norm. In comparison, notations from ``modern mathematics'' as used in everyday practice are substandard, hampering symbolic reasoning and thus making it unpopular. 
   \par The cause of this stagnation is largely historical.  When introducing so-called ``modern mathematics'', forgetting its roots caused serious educational mistakes, denounced in rather strong terms by Arnold~\cite{Arno97otma}.  In a severe overreaction, the view of mathematics as a method was sacrificed in favor of mathematics as a bag of tricks and attempting to elicit motivation by so-called ``real-life'' examples no more realistic than the {\em farmer-sells-potatoes}-type problems in grade school --- and in PISA tests! The well-proven structure {\em definition-examples-theorems} was frowned upon, and mathematical exposition had to become a ``narrative''. 
   \par As a result, classics like Rudin's {\em Principles of Mathematical Analysis}~\cite{Rudi64poma} are, as Krantz observes, ``{\em often no longer suitable, or appear to be inaccessible, to the present crop of students}'' \cite{Kran05raaf}. Here the blame does not fall on the students.
   \par Narratives lack the punctuation provided by headings like ``Definition'' and ``Theorem'', which help novices to distinguish between, say, statements that can be deduced from earlier ones and statements introducing new elements.
   \par If {\em definition-examples-theorems} expositions often deserve criticism, it is not for the usual reasons (take your pick), but because definitions are usually presented as ``given'', or as arbitrary points of departure for a game of logic.  In fact, definitions are the result of {\em design decisions}.  They also determine the flavor of the theorems (and proofs) derived from them.  Hence it is crucial for understanding that these decisions are explained and justified.
   \par In mathematics texts, this is all too rarely done.  One of the few exceptions is Halmos's {\em Naive Set Theory}~\cite{Halm60nast} which, if only for this reason (yet also for other reasons!), should be required reading for all beginning students --- and many mathematicians as well.  Halmos not only explains the design decisions and their shortcomings for most conventions, but also does not shrink back from calling some poor practices ``unacceptable but generally accepted''.  Quine~\cite{Quin69stil} even designates lesser offenses as ``glaring perversity'', which seems an apt characterization of mathematicians acting against better judgment.
   \par  Indeed, perceived ``general acceptance'' is often taken as a licence to perpetuate junk conventions.  Users of inept designs typically defend them by feigning confusion at proper alternatives, calling them ``nonstandard'' even if they have been around for a long time and are routinely used by plenty of other authors. 
   \par If  an engineer is sloppy, his design may fail, even catastrophically. Mathematicians often condone sloppiness, even if it sets bad examples and abuses confidence.  Discerning students will be dissatisfied by the discrepancy between the reputation of mathematics as being precise and actual practice.  Others may even get confused if insecure teachers insist on ``doing things as in the book''. 
   \par Yet, the engineering literature is not blameless either.  Years ago Lee and Varaiya~\cite{LeeV00sisy} corrected many inept mathematical practices in signal processing.
   \par Playing down such issues as ``just a matter of notation'' is misleading. Poor notation prevents the shape of expressions from giving guidance in reasoning.  It also reflects poor understanding, according to Boileau's aphorism ``{\em Ce que l'on con\c{c}oit bien s'\'enonce clairement -- Et les mots pour le dire arrivent ais\'ement}''.  If authors misunderstand their own definitions, what about their students?
   \par ``If it ain't broke, don't fix it'' is another engineering maxim.  Yet, as we shall see, even basic concepts that worked fine 50 years ago somehow got ``broke''.
   \par This paper addresses the issue in the title by presenting a design view on various concepts from the literature, but it is {\em not} some linear, complete proposal.
   \par Often references include page numbers to make them truly useful for the reader.  For brevity, co-authors are omitted when mentioning names in the text.
%
%
\section{Case study A -- Relations: two logically equivalent definitions}
%
%
\subsection{Simple and safe formulations}
%
\noindent The simplest ``modern'' definition of a {\em relation} is typical in older texts such as Bourbaki~\cite[p.~71]{Bour54then},
Suppes~\cite[p.~57]{Supp72axst}, Tarski~\cite[p.~3]{Tars87stwv}, but only in a few current books, such as Jech~\cite[p.~10]{Jech03stth}, Scheinerman~\cite[p.~73]{Sche12madi} and Zakon~\cite[p.~8]{Zako04maan}.
\begin{defin}[Relation]\label{rela}
A relation is a set of ordered pairs.\em
\\Equivalently, in symbols ~\cite{Bour54then,Supp72axst}: $R\;\isrel \siff \forall z \ab z \membr R \simp \exists x \ab \exists y \ab z = (x, y)$\hh.
\end{defin}
Taking {\em set} and {\em ordered pair} colloquially, and with `nonmathematical' examples, the word statement of \dref{rela} is even accessible at grade school level.
   \par In this paper, when saying just ``pair'', we always mean ``ordered pair''.
   \par Some notational design issues arise here. First, an ordered pair is commonly written $(x, y)$. Some authors use $\langle x, y\rangle$, a waste of symbols.  In fact, one can even write $x, y$ and reserve parentheses for emphasis or disambiguation, which also covers $n$-tuples like $(x, y, z)$ and trees like $((x, y), z)$.  Identifying $(x, y, z)$ with $((x, y), z)$ as in Bourbaki~\cite[p.~70]{Bour54then} is clearly a bad design decision. 
   \par Second, the literature diverges about writing $(x, y)$ or $(y, x)$ and $x\,R\,y$ or $y\,R\,x$. Quine~\cite[p.~24]{Quin69stil} offers many good reasons for following Peano and G\"odel in using the {\em natural order} from spoken language, writing ``$a$ is the father of $b$'' as $a\,F\,b$, and ``$a$ is smaller than $b$'' as $a < b$.  Similar reasons would favor writing, for instance, ``velocity versus time'' as $(v, t)$.  However, mathematicians used to writing the ``independent variable'' first might feel disoriented---unlike novices! Tradition can be reconciled with reason by writing $(x, y) \membr R$ iff $y\,R\,x$.  For human engineering reasons, such conventions should be stated conspicuously.
\paragraph{Intermezzo: Quine, the angry notational engineer} In an uncharacteristic diatribe of nearly two pages~\cite[p.~24--26]{Quin69stil}, Quine deplores the ``sorry business'' and ``glaring perversity'' of ill-designed notations.  He concludes (a) ``{\em I have given much space to a logically trivial point of convention because in practice it is so vexatious.}''; (b) ``{\em [Whoever] switched a seemingly minor point of usage out of willfulness or carelessness cannot have suspected what a burden he created.}''.
\spar These remarks reflect typical engineering concerns. Indeed, (a) reminds us that avoiding flaws during the design phase is easy compared with repair afterwards (one line versus lengthy arguments), and (b) emphasizes the importance of taking into account the interests of the users, namely, the future generations.
\paragraph{Auxiliary notions} Most authors using \dref{rela} add \drefs{dora}, \ref{XtoY} and \ref{cmpcnv}.
\begin{defin}[Domain, range]\label{dora}{~}
\\The {\em domain} of a relation $R$ is the set of first members of the pairs in $R$.  The {\em range} of a relation $R$ is the set of second members of the pairs in $R$.\em
\\Notation: the literature mentions various self-explanatory notations such as $\dom R$ or $\mathrm{Dom}(R)$ for the domain of $R$ and $\ran R$ or $\mathrm{Ran}(R)$ for the range of $R$.  
\end{defin}
\begin{defin}[Relation from \boldmath$X$ to $Y$\unboldmath]\label{XtoY}
A {\em relation from $X$ to $Y$}is a relation whose domain is a subset of $X$ and whose range is a subset of $Y$.\em
\\Equivalently, in symbols: $R\,\isrel\;(X, Y) \siff R\;\isrel \sand \dom R \sbset X \sand \ran R \sbset Y$\hh. 
\\Notation: Writing $R : X \lrar Y$ or $R : X \leftrightarrow Y$ introduces a relation from $X$ to $Y$. 
\end{defin}
{\em Aside}\hh As in Meyer \cite[p.~23]{Meye91itpl}, $X \lrar Y$ denotes the set of relations from $X$ to~$Y$.
Also, the symbol~:~clearly distinguishes bindings from statements.  A {\em binding} $i : S$, read ``$i$ \underline{in} $S$'', {\em introduces} an identifier $i$ for an object in a set $S$, whereas  $i \membr S$, read ``$i$ \underline{is in} $S$'' (or similar) is a {\em statement} in which $i$ is {\em used}.  Proper symbolism passes the {\em prose test}: transliterating formulas into words must yield sentences with correct grammar. The RHS of $S \sbset T \siff \evry x \bin S \ab x \membr T$ is read ``for all $x$ \underline{in} $S$, $x$ \underline{is in} $T$'' or ``every $x$ \underline{in} $S$ \underline{is in} $T$''.  Lamport \cite[p.~289]{Lamp03spec} notes that the common forms $\{x \in S \ba p\}$ and $\{e \ba x \in T\}$, where $p$ is a boolean expression and $e$ is any expression, are ambiguous if $p$ is $x \membr T$ and $e$ is $x \membr S$.  Writing $\{x \bin S \ba p\}$ and $\{e \ba x \bin T\}$ yields $\{x \bin S \ba x \membr T\} = S \insec T$ and $\{x \membr S \ba x \bin T\} \sbset \bool$.    \par Finally, $S$ in $i : S$ is called a {\em type} and expresses a range for $i$, not an attribute of $i$. Hence the $X$ and~$Y$ figuring in $R : X \lrar Y$ are attributes of the type $X \lrar Y$, not of $R$.  This concludes the aside.
   \par The following formulation is independent of the $(x, y)$ versus $(y, x)$ issue.
\begin{defin}[Composition, converse]\label{cmpcnv} 
The {\em composition} $S \cmps R$ of relations $S$ and $R$ is the relation such that $z(S \cmps R)x \equiv \exists y . z\,S\,y \wedge y\,R\,x$.
\\The {\em converse} $R\cnv$ of a relation $R$ is the relation defined by $x\,R\cnv\,y \equiv y\,R\,x$.
\end{defin}
Composition is called {\em relative product} by Suppes~\cite[p.~63]{Supp72axst} and Tarski~\cite[p.~3]{Tars87stwv}, and {\em resultant} by Quine~\cite[p.~22]{Quin69stil}.  Suppes writes $S / R$, the others $S | R$.
%
\subsection{Entering murky waters: misconceptions, unsoundness and poor judgement}
%
\noindent All books in our {\em introduction to proof} sample \cite[p.~172]{Bloc11praf}, \cite[p.~101]{Daep03rwap} , \cite[p.~93]{Daep11rwap}, \cite[p.~155]{Garn10dima}, \cite[p.~86]{Gete12imsp}, \cite[p.~51]{Good06dima}, \cite[p.~267]{Grie93logi}, \cite[p.~192]{Hamm09bopr}, \cite[p.~176]{Robe10impt}, \cite[p.~135]{Smit10atta}, \cite[p.~171]{Vell09htpi} except \cite[p.~73]{Sche12madi} combine \drefs{rela} and \ref{XtoY} as follows. 
\begin{defin}[Relation from \boldmath$X$ to $Y$\unboldmath]\label{relcomb}{~} A {\em relation from $X$ to $Y$} is a subset of $X \cart Y$. \em Equivalently, in symbols: $R\;\isrel'\;(X, Y) \siff R \sbset X \cart Y$\hh.
\end{defin}
This forces defining {\em relation} backwards, as a relation from $X$ to $Y$ for some $X, Y$. Also, using Cartesian products at this stage may be an educational burden~\cite{Shua75dima}.  
   \par A more serious problem is that disregarding {\em separation of concerns} carries a heavy price in understanding, apparently even for the authors.  Most fail to recognize that \drefs{XtoY} and \ref{relcomb} are equivalent: $R\;\isrel\;(X, Y) \siff R\;\isrel'\;(X, Y)$.
   \par Indeed, many textbooks strongly suggest that \dref{relcomb} somehow ``glues'' $X$ and $Y$ to the relation, and that one cannot even define just {\em relation} without adding ``{\em from $X$ and $Y$}''. Some texts remain vague here, but the litmus test for one's understanding of a mathematical concept is the view on {\em equality}, in this case: when is a relation $R$ from $X$ to $Y$ equal to a relation $S$ from $U$ to $V$?
   \par The answer was evident 50 years ago: with relations defined as sets, $R = S$ iff both contain the same elements, regardless of $X$, $Y$, $U$, $V$.  Defining equality anew often causes unsoundness, typically by stating that $R = S$ also requires $X = U$ and $Y = V$ \cite[p.~179]{Robe10impt}.  Indeed, if $R \sbset X \cart Y$, $X \neq \mpty$ and $Y \subset V$, then $R \sbset X \cart V$ and $Y \neq V$, so $R = R$ would require $Y = V$, a contradiction. 
   \par Many texts \cite[p.~180]{Robe10impt}, \cite[p.~141]{Smit10atta}, \cite[p.~236]{Vell09htpi}, \cite[p.~94]{Wall12bgdm} define $S \cmps R$ for a relation $R$ from $X$ to $Y$ and a relation $S$ from $U$ to $V$ only for the case $U = Y$.  Since $R \sbset X \cart Y \sbset X \cart (Y \union U)$ and $S \sbset U \cart V \sbset (Y \union U) \cart V$, this is not restrictive, unless one accepts the aforementioned unsound view on equality.
   \par {\em Aside}\hh As a restrictive variant of a {\em relation from $X$ to $Y$}, Bourbaki~\cite[p.~72]{Bour54then} defines a {\em correspondence from $X$ to $Y$} as a triple $(R, X, Y)$ where $R \sbset X \cart Y$, and restricts composition of $(R, X, Y)$ and $(S, U, V)$ to the case $U = Y$.
%
%
\section{Case study B -- Functions}
%
%
\noindent Poor design decisions for relations reappear for functions with a vengeance.
This is especially unfortunate since, as Herstein puts it, {\em Without exaggeration this {\em [namely, a {\em mapping} or {\em function}]} is probably the single most important and universal notion that runs through all of mathematics.}~\cite{Hers64alge}.
   \par The ``modern'' definition of {\em function} was issue-free 50 years ago, and still is in analysis/calculus books~\cite{Apos67calc,Bart64tera,Flet66maan,Kolm70inra,Kran05raaf,Lars09calc,Royd68rean,Rudi64poma,Stew10calc,Sull14calc,Zako04maan}, but unsoundness appears since 2005 in {\em transition to proof} texts \cite{Bloc11praf,Char12mapr,Daep03rwap,Daep11rwap,Garn10dima,Gete12imsp,Hamm09bopr,Robe10impt,Smit10atta}.
   \par In passing, we mention some harmful myths that can be read between the lines in textbooks but surface explicitly in oral and written conversations. Myth \#0 (a meta-myth actually) holds that divergences in definitions just reflect different needs in various disciplines \cite{Shua75dima}.  However, our samples come from algebra, analysis/calculus, discrete math, logic, set theory, and reveal that nearly {\em all} of them use the {\em same} function concept, differing only in the care devoted to design and formulation. Myth \#0 is harmful in trying to divert closer scrutiny. 
   \par Recurrent points of interests are: (i)~defining {\em function}, (ii)~function equality, (iii)~function from $X$ to $Y$, (iv)~onto-ness, (v)~function composition and inverse. 
%
\subsection{Once again: simple and safe formulations}
%
\pbfni{(i)} A relation $R$ is called {\em functional}~\cite{Bour54then,Meye91itpl} iff no two pairs in $R$ have the same first member.   Hence the following phrasings are equivalent; the choice depends on whether {\em relations} are skipped, as in analysis/calculus texts, or defined first. 
\begin{defin}[Function]\label{AB} 
{\em (A) Apostol~\cite[p.~53]{Apos67calc}:} A {\em function} $f$ is a set of ordered pairs $(x, y)$ no two of which have the same first member.
\\{\em (B) Bourbaki~\cite[p.~77]{Bour54then}:} A {\em function} is a functional relation.
\end{defin}
\dref{AB} is also found in Dasgupta \cite[p.~10]{Dasg14sett}, Flett~\cite[p.~4]{Flet66maan}, Jech \cite[p.~11]{Jech03stth}, Mendelson \cite[p.~6]{Mend87itml}, Scheinerman \cite[p.~167]{Sche12madi}, Suppes \cite[p.~86]{Supp72axst}, Tarski \cite[p.~3]{Tars87stwv}, Zakon~\cite[p.~10]{Zako04maan}. Functionality justifies writing $y = f(x)$ iff $(x, y) \membr f$.  As in \cite[p.~1]{Macl71cfwm}, one may write $f\,x$ instead of $f(x)$ when no ambiguity results. Authors using \dref{AB} introduce the domain and range as in \dref{dora}.
   \parbf{(ii)} This results in the following theorem, quoted from Apostol \cite[p.~54]{Apos67calc}.
\begin{therm}[Equality]\label{apotherm}
Functions $f$ and $g$ are equal iff \hh {\nrm(a)} $f$ and $g$ have the same domain, and {\nrm(b)} $f(x) = g(x)$ for every $x$ in the domain of $f$.
\end{therm}
   \parbf{(iii)} Authors starting from \dref{AB}, including Apostol~\cite[p.~578]{Apos67calc}, Dasgupta \cite[p.~10]{Dasg14sett}, Flett~\cite[p.~5]{Flet66maan}, Jech \cite[p.~11]{Jech03stth}, Scheinerman \cite[p.~169]{Sche12madi} and Zakon \cite[p.~10]{Zako04maan}, use the following common notions for classifying functions.    
\begin{defin}[Function from \boldmath$X$ to $Y$\unboldmath]\label{fXtY} A {\em function from $X$ (in)to $Y$} is a function with domain~$X$ and range included in~$Y$.\em
\\Notation: writing $f : X \arr Y$ introduces a function $f$ from $X$ to $Y$.
\end{defin}
This is the ISO standard~\cite[p.~15]{ISOs09quun}, where a {\em function} is defined in broader terms, mentioned later.  Divergent views on $f \bin X \arr Y$ are {\em nonstandard}. 
   \par  It is convenient reading $X \arr Y$ as the set of functions from $X$ to $Y$ and $X \parr Y$ as the set of functional relations from $X$ to $Y$ \cite[pp.~25--26]{Meye91itpl}.  Such types serve as {\em partial specifications} for functions. Tighter types are defined later.  
   \parbf{(iv)} Next, we consider {\em onto-ness} as defined by authors using \drefs{AB} and \ref{fXtY}, for instance, Flett \cite[p.~5]{Flet66maan}, Jech \cite[p.~11]{Jech03stth}, Mendelson \cite[p.~6]{Mend87itml}, Scheinerman \cite[p.~172]{Sche12madi}, Tarski \cite[p.~3]{Tars87stwv} and Zakon \cite[p.~11]{Zako04maan}.
\begin{defin}[Onto \boldmath$Y$\unboldmath]\label{ontoY} For any set $Y$, a function is said to be {\em onto $Y$}, or {\em surjective on $Y$}, iff its range is $Y$.
\end{defin}
Note that {\em onto} is a preposition, and appears as such in \dref{ontoY}, which is also used by many authors (listed later) who do {\em not} start from \dref{AB}.
   \par The dual notion of ``$f$ is {\em onto $Y$}'' ($\tran f = Y$) is ``$f$ is {\em total on $X$} ($\tdom f = X$).  The dual notion of ``$f$ is {\em into $Y$}'' ($\tran f \sbset Y$) is ``$f$ is {\em partial on $X$}'' ($\tdom f \sbset X$).  A~function $f : X \parr Y$ is often called a {\em partial function}, but is $\{(0, 1), (2, 3)\}$ a partial function? Grammatically correct is: a {\em function from part of $X$ to $Y$}.
   \parbf{(v)} For composition, we mention two equivalent formulations. Formulation (A) skips relations, as in Apostol~\cite[p.~140]{Apos67calc}, Flett~\cite[p.~11]{Flet66maan}, Mendelson~\cite[p.~7]{Mend87itml} and many others, listed later.  Formulation (B) is based on relations.
\begin{defin}[Function composition \boldmath$g \cmps f$\unboldmath]\label{compos} {\em (for {\em any} functions $g$ and $f$)}
\\{\em(A)} $g \cmps f$ is the function whose domain consists of all $x$ in $\tdom f$ that satisfy $f(x) \membr \tdom g$ and whose value $(g \cmps f)(x)$ for arbitrary $x$ in that domain is $g(f(x))$.
\\{\em(B)} \em $g \cmps f$ follows \dref{cmpcnv} assuming natural order as in Quine~\cite[p.~24]{Quin69stil}; otherwise $f$ and $g$ must be swapped, e.g., $g \cmps f = f / g$ in Suppes~\cite[p.~87]{Supp72axst}, $g \cmps f = f | g$ in Tarski~\cite[p.~3]{Tars87stwv}.  Proof obligation: showing that $g \cmps f$ is functional.
\end{defin}
%
\subsection{Wearing the ice thin: convoluted formulations}
%
\pbfni{(i--iii)} \drefs{AB} and \ref{fXtY} are sometimes crammed together, starting as early as Halmos \cite[p.~30]{Halm60nast} and Herstein \cite[p.~10]{Hers64alge}, and more often in current texts including Krantz \cite[p.~20]{Kran05raaf}, Velleman \cite[p.~226]{Vell09htpi} and others, listed later.
\begin{defin}[Function from \boldmath$X$ to $Y$\unboldmath]\label{C} Let $X$ and $Y$ be sets.
\par {\em (a)} A {\em function $f$ from $X$ to $Y$}, written $f \bin X \arr Y$, is a relation $f \sbset X \cart Y$ satisfying the property that for each $x$ in $X$ the relation $f$ contains exactly one ordered pair of the form $(x, y)$.
\par {\em (b)} The set $X$ is called the {\em domain of $f$}.
\end{defin}
As a warning against uncritical copying, we reproduced the widespread but unacceptable phrasing, which suggests that the function is written $f \bin X \arr Y$ (in fact, the function is written just $f$), and that $f \sbset X \cart Y$ is a relation (in fact, it is a statement about the relation $f$).  Proper  phrasings are evident.
%
%
   \par All authors using \dref{C} overlook that part (b) requires proving that $X$ is fully determined by $f$ as defined in (a).  This is easy; the result is $X = \tdom f$.
   \par More importantly, users of \dref{C} fail to realize its logical equivalence to \dref{AB} (proof: exercise), including the standard meaning of $f : X \arr Y$.  Still, the different formulation has a huge impact on clarity.  By disregarding separation of concerns, \dref{C} has given rise to Myth \#1, which holds that one cannot define {\em function} by itself, but only {\em function from $X$ to $Y$}.
   \par The common pitfalls are exposed by the litmus test: equality. Surprisingly few present-day texts using \dref{C} mention \tref{apotherm}, found only in Daepp~\cite[p.~152]{Daep03rwap}, Gerstein \cite[p.~113]{Gete12imsp} and Smith et al.~\cite[p.~189]{Smit10atta}.  Instead, many define equality anew, all too often unsoundly, as demonstrated later.
   \parbf{(iv)} ``Classical'' authors using \dref{C} (or similar), including Bartle \cite[p.~13]{Bart64tera}, Halmos \cite[p.~31]{Halm60nast}, Herstein \cite[p.~12]{Hers64alge}, Kolmogorov \cite[p.~5]{Kolm70inra}, say that $f$ is {\em onto~$Y$} iff $\tran f = Y$, as in \dref{ontoY}, using ``onto'' as a {\em preposition}.  Some of the few ``modern'' users of \dref{C} who write ``onto $Y$'' are Gerstein \cite[p.~118]{Gete12imsp} and Smith \cite[pp.~xvii, 205]{Smit10atta}, but their formulation lacks generality. 
   \parbf{(v)}  Most ``classical'' authors using \dref{C} (or similar), including  Bartle~\cite[p.~40]{Halm60nast} and Halmos~\cite[p.~40]{Halm60nast}, define $g \cmps f$ for arbitrary functions $f$ and $g$, as in \dref{compos}.  Only a few classical texts \cite[p.~13]{Hers64alge}\cite[p.~9]{Royd68rean} restrict coverage of $g \cmps f$ for functions $f : X \arr Y$ and $g : U \arr V$ to the special case $U = Y$.
   \par Most ``modern'' texts succumb to this restriction, including Bloch \cite[p.~146]{Bloc11praf},  Roberts \cite[p.~226]{Robe10impt},  Scheinerman \cite[p.~183]{Sche12madi}, Smith \cite[p.~197]{Smit10atta}, Velleman \cite[p.~231]{Vell09htpi}.  An intermediate form requiring $\tran f \sbset \tdom g$ appears in Daepp (2003) \cite[p.~175]{Daep03rwap}, Daepp (2011) \cite[p.~167]{Daep11rwap}, Jech \cite[p.~11]{Jech03stth} and Krantz \cite[p.~22]{Kran05raaf}.  The general form appears in Larson \cite[p.~25]{Lars09calc} and Stewart \cite[p.~40]{Stew10calc} which, not surpringly, are calculus texts, since restricted composition is impractical. 
%
\subsection{Falling through the ice: common yet unsound additions to \dref{C}}
%
\pbfni{(i--iii)} \dref{C} is the sound part of definitions in {\em transition to proof} texts, but Bloch \cite[p.~131]{Bloc11praf}, Chartrand \cite[p.~216]{Char12mapr}, Daepp \cite[p.~147]{Daep03rwap}\cite[p.~143]{Daep11rwap}, Garnier \cite[p.~224]{Garn10dima} Gersting \cite[p.~383]{Gete12imsp},  Hammack \cite[p.~195]{Hamm09bopr}, Roberts \cite[p.~220]{Robe10impt}, Smith \cite[p.~185]{Smit10atta},  Gilbert \cite[p.13]{Gilb08elma} and Wallis \cite[p.~106]{Wall12bgdm}, add
\begin{defin}[Codomain]\label{cod} {\em \dref{C}}{\em(c)} $Y$ is called the {\em codomain of $f$}.
\end{defin}
However, just like \dref{C}(b), adding 10(c) entails a proof obligation.  Recognizing this clearly reveals a logical contradiction.  Indeed, the definiendum is {\em codomain of $f$}, the definiens is $Y$, but $Y$ is not uniquely determined by $f$. As for relations, letting $f \sbset X \cart Y \subset X \cart Y'$ reveals a contradiction.   
   \par Still, some authors uphold Myth \#2: writing $f : X \arr Y$ makes $Y$ an attribute of $f$ by specifying $f \sbset X \cart Y$.  Yet, all this says about $Y$ is $\tran f \sbset Y$.
   \par Myth \#3 maintains that \dref{C} contains ambiguities allowing multiple views, making logical contradictions just a matter of interpretation.  Yet, in \dref{C}(a), the definiendum and the definiens are clear (except as written in \cite[p.~131]{Bloc11praf}), using the unambiguous concepts {\em subset} and {\em Cartesian product}.
   \par Also, insofar as \dref{cod} makes the perceptive reader wonder if the authors really mean ``codomain {\em of $f$}'', further context indicates they mostly do. 
   \par Again equality is most revealing.  Apart from three exceptions mentioned, all {\em transition to proof} texts using \dref{cod} overlook \tref{apotherm} and define equality anew.  Roberts \cite[p.~223]{Robe10impt} avoids conflict with \dref{C} by using the statement of \tref{apotherm}.  Others, e.g., Bloch \cite[p.~136]{Bloc11praf}, Garnier \cite[p.~224]{Garn10dima}, Hammack \cite[p.~198]{Hamm09bopr} extend this statement with $\codt f = \codt g$, indirectly contradicting \dref{C}.
%
%
Interestingly, Smith \cite[p.~189]{Smit10atta} explicitly states that function equality does {\em not} require equal codomains! 
   \parbf{(iv)} Unsoundness also results from improper use of ``{\em onto}'' as an adjective, as in Bloch \cite[p.~155]{Bloc11praf}, Daepp \cite[p.~163]{Daep03rwap}\cite[p.~157]{Daep11rwap}, Gries~\cite[282]{Grie93logi}, Hammack \cite[p.~199]{Hamm09bopr}, Krantz \cite[p.~22]{Kran05raaf},
Roberts \cite[p.~231]{Robe10impt} and Velleman \cite[p.~236]{Vell09htpi}.
\begin{defin}[Onto]\label{onto} A function $f : X \arr Y$ is {\em onto} {\em (surjective)} iff $\tran f = Y$.
\end{defin}
With \dref{onto}, the same function can be both onto and not onto depending on whether or not the set $Y$ appearing in $f : X \arr Y$ happens to be $\tran f$.
   \parbf{(v)} All texts adding \dref{cod} require for $g \cmps f$ that $\tdom g = \tcod f$ and for the inverse $f^-$ that $\tran f = \tcod f$. This is impractical for applications, e.g., in calculus.
%
\subsection{Design considerations, variant concepts, and evaluation}
%
\noindent Many textbooks and countless blogs use the term {\em codomain}, typically in the unsound way described.  Even though most authors using \dref{C}/\ref{cod} seem to sense the problems with squeezing in codomains, they somehow feel obliged to try. Of course, they needn't! Clearly, a proper account for {\em codomain} is overdue. 
   \par In view of the earlier analysis, the term {\em codomain} is best (i) simply discarded, or (ii) used as the symmetric counterpart of {\em domain}, thus far called {\em range}, or (iii) recognized as an attribute of a type like $X \arr Y$ or $X \parr Y$, not of a function. 
   \par A quite different approach is defining a variant of the function concept such that the set $Y$ in $f : X \arr Y$ is truly part of $f$, safely called {\em the codomain of~$f$}.
   \par For instance, Bourbaki initially defines a function \cite[p.~76]{Bour54then} as a triple $(F, A, B)$ where $F$ is a functional relation with domain $A$ and range included in~$B$.  Such a triple is subsequently called an {\em application from $A$ into $B$} \cite[p.~76]{Bour54then}, which avoids confusion with using {\em function} for a functional relation \cite[p.~77]{Bour54then}.  The term {\em codomain} is not mentioned in \cite{Bour54then} --- so let's not blame Bourbaki! 
   \par The engineering question is: what purpose might such a variant serve?
   \par For the sake of generality, this issue is best disassociated from the set of pairs view, whose predominance in the discussion thus far just reflects random literature samples from diverse areas of mathematics.  Many authors, including Lang~\cite[p.~38]{Lang83unan}, Lee~\cite[p.~48]{LeeV00sisy},  Royden~\cite[p.~8]{Royd68rean} and Spivey~\cite[p.~29]{Spiv89thzn}, note that a set of pairs is really a \emph{representation} of a function, called its \emph{graph}. Hence let's broaden the {\em representational} \dref{AB} to a {\em conceptual} one, inspired by the Goursat/Courant style, but generalized to arbitrary domains and properly distinguishing $f$ from $f(x)$.  For instance, the ISO standard~\cite[p.~15]{ISOs09quun} just says that a {\em function} $f$ assigns to each $x$ in its domain a unique value $f(x)$.  Here ``assigns'' can be made precise by an {\em assertion} of the form $A(x, f(x))$, called {\em relation} by Bourbaki \cite[p.~47]{Bour54then}, but not to be confused with a set of pairs. 
   \par {\em Equality} is pivotal for mathematical objects. Distinguishing functions on the basis of possible assignments outside their domain would be an useless complication, hence these values are best declared irrelevant for a function.  
   \par (A) The minimalist, ``no frills'' design reflecting this view is the following.
\begin{defin}[The {\em function} concept: minimalist design]\label{funcon}{~}
\par {\em i.} A {\em function} $f$ is an object fully specified by {\em(a)}~a set $\tdom f$, called {\em the domain of $f$}, and {\em(b)}~for each $x$ in $\tdom f$ a unique {\em value}, written $f(x)$ or $f\,x$.
\par {\em ii.} The stipulation ``{\em fully specified}'' means that $f = g$ if {\em(a)}~$\tdom f = \tdom g$ and {\em(b)}~$f(x) = g(x)$ for all $x$ in $\tdom f$. {\em(Note: ``only if'' by Leibniz's principle \cite{Grie93logi})}
\par {\em iii.} A {\em function $f$ from $X$ to $Y$} is a function such that {\em(a)}~$\tdom f = X$ and {\em(b)}~$f(x) \membr Y$ for all $x$ in $X$.  Such a function is introduced by writing $f : X \arr Y$.
\end{defin}
\dref{funcon}.iii simply follows the ISO standard: $\tdom f = X$ and $\tran f \sbset Y$, where $\tran f = \{f(x) \ba x \bin \tdom f\}$. Types of the form $X \arr Y$ are partial specifications.  An illustration is defining $\mathrm{sqrt} \bin \real_{\geq 0} \arr \real_{\geq 0}$ with $(\mathrm{sqrt}\,x)^2 = x$. {\em Composition} is unrestricted: $\dom (g \cmps f) = \{x \bin \tdom f \ba f(x) \membr \tdom g\}$ and $(g \cmps f)\,x = g(f\,x)$ as usual.
   \par Composition also supports {\em specification by proxy}: specifying $g \bin Y \arr Z$ via $f \bin X \onto Y$ ($\onto$ indicating onto) and $h \bin X \arr Z$ by $g(f\,x) = h\,x$. {\em Well-definedness} (functionality) of $g$ requires $h(x) = h(x')$ whenever $f(x) = f(x')$.
   \par An example is defining the inverse: let $g := f^-$, $Y := \ran f$, $Z := X$ and  $h := \ide_X$. Well-definedness of $f^-$ amounts to $f$ being 1-1.  In the general scheme, if $f$ is 1-1, then $g$ is well-defined and $g \cmps f = h$ is equivalent to $g = h \cmps f^-$.  Pattern matching is an instance: compare $g(\mathrm{cons}(a, x)) = h(a, x)$ and $g\,s = h(\mathrm{cons}^- s)$.
   \par The set of pairs view will remain a useful analogy, e.g., in defining $f^-$ as the object represented by the inverse relation, which is functional iff $f$ is 1-1.
   \par (B) A typical non-minimalist design variant of \dref{funcon} would add: 
 i.(c)~a~set $\tcod f$, called {\em the codomain of $f$};\hh\ 
 ii.(c) $\tcod f = \tcod g$;\hh\ iii.(c) $\tcod f = Y$.
   \par Logically, composition and inverses {\em could} still be defined without restriction (exercise).  However, authors using codomains do restrict $g \cmps f$ by $\tdom g = \tcod f$ and define inverses for 1-1 functions only if the latter are ``onto'' [their codomain]. 
   \par So the issue boils down to: what are the costs and the merits of codomains?  
   \par Shuard~\cite{Shua75dima} published perhaps the only paper evaluating the function-with-codomain variant.  Her single (!) argument {\em in favor} is using {\em onto} as an adjective. Yet, the ability to say that ``$f$ is onto $Y$ but not onto $Z$'' is more selective.
   \par Shuard's argument {\em against} is more solid: simplicity.  She states: ``{\em Flett's definition wins hands down as simplicity in analysis is concerned}''. Her next statement, ``{\em In algebra, however, it is more convenient to start by mentioning the codomain of a function}'' (resembling Myth \#0), is left unsubstantiated.  Excellent algebra texts such as Herstein~\cite{Hers64alge} do fine with the standard variant.  
   \par In Shuard's paper and all other sources consulted, suggestions that attaching a codomain {\em might} be convenient turns out to be fallacious, typically overlooking that the standard view regarding $f : X \arr Y$ already implies $\tran f \sbset Y$.  This view covers all sensible purposes of $f : X \arr Y$, ``mentioning'' $Y$ included, and without making $Y$ a function attribute.  Only for topology further study is needed to determine whether viewing $Y$ as an attribute of a function $f : X \arr Y$ is just some tradition based on  similar oversights or has genuine advantages.
   \par Still, what's the harm in a function-with-codomain, beside complexity? 
   \par From a conceptual and practical viewpoint, burdening a function with a codomain affects all other definitions, complicates equality and impoverishes the function algebra for inverses, composition, merge, override and so on \cite{Meye91itpl}.  Shuard~\cite[p.~10]{Shua75dima} notes that the only analysis book that she knows to use codomains~{\cite{Spre70eora} gets into trouble by defining the inverse of a 1-1 function $f$ to have domain $\tran f$ as usual and ignoring that codomain users require surjectivity.  She adds that ``{\em the distinction between $f : A \arr \real$ and $f_1 : A \arr f(A)$ is so tedious that it is clearly better forgotten at this stage}''.  One might say: ``{\em clearly better avoided from the start}'', matching Smith's view on equality \cite[p.~189]{Smit10atta}. 
   \par {\em Aside: programming versus mathematics} \hh Types and signatures of ``functions'' in programming~\cite[Section 3.8]{Garn10dima} and some proof assistants typically are unique attributes by design.  They are easier to implement than general symbolic computation of, say, $\dom (g \cmps f)$, but remain rather crude approximations of types as partial specifications following the ISO standard.  Indeed, the generality provided free of charge by the minimalist/standard view (illustrated  for $g \cmps f$ and $f^-$) is common fare in {\em mainstream mathematics}, by which we mean: the mathematics routinely applied by the large majority of users, ranging from mathematicians active in analysis/calculus, linear algebra, discrete math etc. to engineers active in signals and systems.  Such users would be justified in dismissing as impractical any definition that infringes on these ``acquired rights''.
   \par Similar considerations in the context of specification languages are found in~\cite{Lamp99slan}.  In the specification language TLA$^+$~\cite[p.~48]{Lamp03spec}, the notation $[X \arr Y]$ has the meaning of $X \arr Y$ as defined by the ISO standard.
   \par Not surprisingly, nearly all calculus/analysis texts avoid codomains and simply proceed from \dref{AB} or \ref{funcon} or equivalent, the most complete picture being presented by Flett~\cite[pp.~4--6]{Flet66maan}.  In this manner, calculus/analysis books succeed in giving a proper account in about one page, without being too terse, and often before page 10.  Functions are too important to postpone their introduction beyond page 100, only to get them bogged down in unsoundness.
   \par {\em Conclusion}\hp2 For mainstream mathematics, attaching a codomain onto a function has no verified merits but increases complexity  and reduces generality.  Hence any definition that accepts such drawbacks entails a heavy obligation of justifying the design, even if it concerns only some niche area.
%
%
\section{Engineering mathematical abstraction}
%
%
\noindent Abstraction is a very useful intellectual tool in science, especially engineering. It allows reasoning about the essentials, without sidetracking by elements causing unnecessary complications. Of course, {\em unnecessary} is the operative word.
   \par Being a tool, abstraction requires engineering. Criteria include not just soundness, but also effectiveness in reasoning. Parnas~\cite{Parn90educ} deplored 25 years ago that ``{\em Those working in theoretical computing science lack an appreciation for the simplicity and elegance of mature mathematics}''.   Little has changed. The context in \cite{Parn90educ} makes clear that {\em mature} refers to mainstream mathematics as characterized earlier. Symbolic reasoning has been kept from reaching maturity by decades of stifling in favor of ``narrative'' prose~\cite{Vell09htpi}, but that is another story.
   \par Engineering mathematical abstraction involves {\em balance}. Abstraction easily degenerates into obfuscation, making it ineffective as a tool for general use.
%
%
\section{Illustration A: functions and the Halmos principle}
%
%
\noindent The function concept is the result of a long evolution.  Its ``modern'' formulation as in \dref{AB} is a one-liner in the most positive sense: a degree of simplicity and clarity that is unlikely to be surpassed. This is typical for a final design. Arguably, non-equivalent variants need a different name to avoid confusion%
\footnote{Halmos~\cite{Halm57bour} is quite amused by Bourbaki's habit of abandoning their ``innovations'' in favor of common terminology. Still, Bourbaki's term {\em application} for the ``triples'' variant \cite[p.~76]{Bour54then} at least avoids confusion with using {\em function} for the common variant \cite[p.~77]{Bour54then}.}.
   \par Yet even when using the set-of-pairs definition as a reference basis, in practice one rarely thinks about functions as sets, and set theory serves only as a handy framework, as noted by Halmos  \cite[p.~31]{Halm60nast}, and illustrated by \dref{funcon}.
   \par Reducing concepts to a set-theoretic representation is conceptually chafing and causes ``freak properties'' or ``accidental facts''  \cite[pp.~25, 45]{Halm60nast}.  Two examples: defining a function as a set (of pairs) yields $(x,y) \membr f$, and defining an ordered pair $(a, b)$ as $\{\{a\}, \{a, b\}\}$ yields $\{a, b\} \membr (a, b)$.  Halmos considers such effects ``{\em a small price to pay for conceptual economy}'' and offers two solutions.
   \par One solution is axiomatization.  Bourbaki~\cite[p.~68]{Bour54then} characterizes pairs by the equality axiom $(x, y) = (x', y') \simp x = x' \sand y = y'$.  \dref{funcon} characterizes functions by the equality axiom $f = g \Leftarrow \tdom f = \tdom g \sand \evry x \bin \tdom f \ab f(x) = g(x)$. In both examples, the converse is trivial by Leibniz's principle~\cite[p.~60]{Grie93logi}.
   \par The other solution consists in using overly concrete definitions only to derive theorems that capture the essence of a concept, and then declaring ``{\em the definition has served its purpose by now and will never be used again}'' \cite[p.~25]{Halm60nast}.  Although many others tacitly follow this approach as well, we name it the {\em Halmos principle} after the author who made it explicit as a design strategy.
   \par According to Halmos, ``{\em The mathematician's choice is between having to remember a few more axioms and forgetting a few accidental facts; the choice is pretty clearly a matter of taste}''.  Still, expressiveness and support for reasoning are more reliable criteria.  The choice also depends on each specific concept.
   \par For instance, distinguishing functions from sets as in \dref{funcon} helps disambiguating common expressions such as $f^n$ and $S^n$.  This is secondary. 
    \par For reasoning, the ``set of  pairs'' view remains a powerful simplifying analogy, as shown for composition and inverses.  It is linked to \dref{funcon} as the {\em graph} of a function: $\tgrph f = \{(x, f\,x) \ba x \bin \tdom f\}$. Clearly, $f = g \simp \tgrph f = \tgrph g$ by Leibniz's principle. By \dref{funcon}.ii, $\tgrph f = \tgrph g \simp f = g$.  Hence $\tgrph$ is invertible and can be omitted (left implicit) when using the ``set of pairs'' view as an analogy.
%
%
\section{Illustration B: Preparing category theory for engineering}
%
%
\subsection{Categories and the Arnold Principle}
%
\noindent The practical potential of category theory has been demonstrated in~\cite{Bird97aofp}.  This book, as well as \cite{Arbi75asaf}, are among the best introductions to the subject.
   \par Still, an engineer would ask: {\em why} use category theory? A non-evasive answer is: to support the algebraic elegance and practical benefits of the point-free style, which has proven useful throughout engineering. For instance, in systems modeling, this style has helped making common notations more general~\cite{Bout03cgf} and free of errors that ``undermine the students' confidence in mathematics''~\cite{LeeV00sisy}.
   \par Why {\em category} theory? Because it is ready-made and seems plug-and-play.
   \par Yet, the literature based on this view reveals grave mismatches. Category theory is about {\em arrows}.  In the {\em strict} (i.e., textbook) variant, each arrow $a$ has a unique {\em source} $\src a$ and {\em target} $\tgt a$. {\em Composition} $b \cmps a$ requires $\src b = \tgt a$.  In this setup, a {\em function} $f$ corresponds to many {\em arrows} (one for each triple $f, A, B$), as noted by Bird \cite[p.~26]{Bird97aofp} and Pierce \cite[p.~2]{Pier91ctcs}.  The target attribute inherits all burdens of the codomain.  Hence, in strict category theory, arrows do not fit {\em relations} and {\em functions} but {\em correspondences} and {\em applications} instead.
   \par Pierce \cite[p.~3]{Pier91ctcs} adds that arrows fit functions in their ``ordinary mathematical meaning'' only in variants of category theory where one cannot always tell for an arrow what its source and target are.  But uniqueness of source and target is precisely why users of strict category theory feel forced to restrict composition!
   \par Such design decisions bring to mind the following familiar parable. 
   \par {\em The Suit}\hp1 A man buys a new suit. When he finds it too tight at some places and too loose at others, the tailor shows him how twisting his arms, turning his feet and bending backwards brings relief. When the customer proudly hobbles along in his new suit, a passerby observes ``How sad, this man must have had a terrible accident''.  Says his companion, ``Yes, but he is lucky to have found a true master tailor who can make such a perfectly fitting suit''.
   \par The moral is that well-designed axiomatization is fitted to existing concepts, not the other way around like a straitjacket.  Typical examplary designs are group, ring and field theory \cite{Hers64alge}: they capture commonalities of existing concepts of quite different nature, without distorting them to fit the axioms.
   \par Arnold~\cite{Arno97otma}, one of the greatest 20th century mathematicians, denounces improper use of abstraction by ``the criminal algebraists-axiomatisators''.  Despite this hyperbole, one can hardly deny his observation that often ``{\em the so-called `axioms' are in fact just (obvious) properties} [of the concepts of interest]''.  This suggests an apt guideline, here called the {\em Arnold Principle}: unless an axiomatization is applied regularly to objects more abstract than those of primary interest, it is obfuscation.  Arnold mentions Whitney's theorem to warn that more abstract objects may not always exist.  Even if they do exist, proper axiomatization unifies them with the objects of primary interest, but never constrains the latter, which would amount to upside-down design.  Inability to adequately capture functions, ``the single most important and universal notion that runs through all of mathematics''~\cite{Hers64alge}, inhibits the practicality of category theory.
   \par Rather than going for extremes such as dismissing either category theory or its flaws, we show how categorical concepts can suit mainstream mathematics.
%
\subsection{Reaping the benefits of category theory without the constraints}
%
\noindent Flawed axiomatization by upside-down design is most effectively avoided by starting from the objects of interest and deriving theorems reflecting abstract properties \cite{Arno97otma,Bout88sacs}.  When appropriate (by Arnold's principle), such theorems can be recast as axioms afterwards.  The proofs then become evidence that the objects of interest do satisfy these axioms.  Thus, no effort is wasted.
   \par Proper design avoids misleading terminology by using {\em relation} and {\em function} only in their ordinary mathematical meaning.  For the concepts captured by strict category theory, safe terms are {\em correspondence} and {\em application} (or alternatives).  Abstractions of correspondences are often called {\em allegories} \cite{Bird97aofp,Frey90caal}.
   \par The relevant benefits of category theory are point-free expression and reasoning.  Still, good design should not cause a rift between point-free and point-wise styles, since practical applicability requires safe and smooth mixing of styles.
   \par One approach along these design guidelines is {\em concrete category theory}, which starts from the primary objects of interest (relations and functions as usual, with the standard view on $f : X \arr Y$), and hence is ``ordinary mathematics''.
   \parbf{(i)} We assume a {\em relation} $R$ is defined conceptually in the style of \dref{funcon} (exercise), with the set of pairs view as an analogy: $y\,R\,x \siff (x, y) \membr \tgrph R$ and $R = S \siff \tgrph R = \tgrph S$.  The domain $\dom R$ and the range $\ran R$, defined via \dref{dora}, are proper attributes of $R$.  Bringing them into the picture often yields sharper results ($=$ instead of $\sbset$). Composition $S \cmps R$ (or $R\;; S$) and converse $R\cnv$ follow \dref{cmpcnv}.  The {\em identity relation $\ide_A$ on $A$} with $\ide_A = \rna{(a, a) \ba a \bin A}$. Predicate calculus yields the usual collection point-free formulas.
   \parbf{(ii)} A relation $R$ is {\em functional} iff $R \cmps R\cnv = \ide_{\tran R}$ (equivalently, $R \cmps R\cnv \sbset \ide_Y$ for $R : X \lrar Y$).  Such a relation is called a {\em function}, typically written $f$, $g$, $\ldots$.
   \par Symbolism can be augmented by diagrams of arrows labeled by relation or function expressions, e.g., $\tlngr{f} \tlngr{g}\;=\;\tlngr{f;g}$.  Optionally, types can be specified by labels at endpoints of arrows and shaping the latter like the type arrows $\tlrr$, ${\arr}$ and $\tparr$, more conveniently drawn as \tdarr. Here are some illustrations:
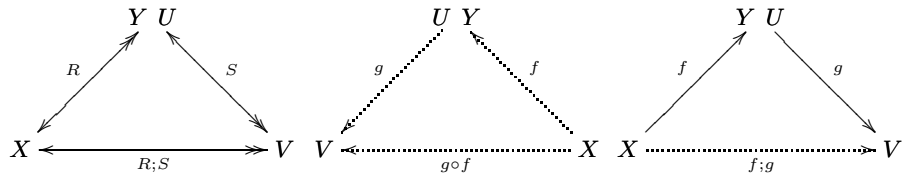
\[\bfig\small
\Atriangle/{<<->}`{<->>}`{<->>}/[Y\hh U`X`V;R`S`R ; S]
\Atriangle(1150,0)/{.>}`{<.}`{<.}/[U\hh Y`V`X;g`f`g \cmps f]
\Atriangle(2300,0)/{<-}`{->}`{.>}/[Y\hh U`X`V;f`g`f ; g]
\efig\]
Recall that types are partial specifications for, not attributes of, the arrows.  Normally, extra structure naturally merges endpoint labels at a common vertex. Then diagrams look exactly as in strict category theory, but without  constraints.  In the following examples, it will become evident that attaching a target or codomain to a function or relation would only provide negative added value.
%
\subsection{Some illustrations of concrete category theory --- selected topics}
%
\noindent{\em Convention}\hh In some contexts, {\em family} is a graphic synonym for {\em function} \cite[p.~77]{Bour54then} \cite[p.~34]{Halm60nast}. An {\em $I$-family} is a function with domain $I$, called {\em index set}.
\paragraph{Product}
The first topic has significant fundamental and practical interest. 
\begin{defin}[Product]\label{prddfn}
The {\em Cartesian product} $\Cart T$ of a family $T$ of sets is the set of functions $f$ such that $\tdom f = \dom T$ and $f(x) \membr T(x)$ for all $x$ in $\tdom f$.
\end{defin}
Products support {\em dependent types} and {\em tolerances on functions} \cite{Bout03cgf} that can be tight, e.g., if $T$ is a family of singletons, declaring $f : \Cart T$ fully specifies $f$.  Hence types as partial specifications can be arbitrarily fine, and $\Cart T \sbset \dom T \arr \Union T$.
   \par Another example is $T := (X, Y)$.  An $n$-{\em tuple} is a function whose domain consists of the first $n$ natural numbers \cite[p.~45]{Halm60nast}\cite[p.~75]{Kran05raaf}\cite[p.~9]{Royd68rean}.  For an $n$-tuple $T$ of sets, one can use infix notation, e.g., $\Cart(X, Y) = X \cart Y$.  This view also resolves the ``terminological friction'' discussed by Halmos \cite[p.~36]{Halm60nast}.
   \par To obtain the usual point-free characterization of $\Cart T$,  we let $I := \dom T$ and define an $I$-family $\pr$ of {\em projection functions} by $\pr_i \membr \Cart T \arr T_i$ and $\pr_i t = t_i$. Note: an illustrative equivalent declaration is $\pr \bin \Cart_{i \bin I} (\Cart T \arr T_i)$.
\begin{therm}\label{prdthm} Let $T$ be an $I$-family of sets, all nonempty {\em (to avoid $\Cart T = \emptyset$)}.  
\pbfni{\em(i)} For any set $S$ and $I$-family $f$ of functions with $f_i \membr S \arr T_i$, there exists a unique $g : S \arr \Cart T$ satisfying $f_i = \pr_i \cmps g$. Specifically, $g = f\trsp$ {\em (transpose)}.
\pbfni{\em(ii)} Let $\gam$ be an $I$-family of functions with $\gam_i \membr C \arr T_i$ and the property that, for any set $S$ and $I$-family $f$ of functions with $f_i \membr S \arr T_i$, there exists a unique $h : S \arr C$ satisfying $f_i = \gam_i \cmps h$. Then there is a bijection between $C$ and $\Cart T$.
\em
\pbfni{Proof}\hh \textbf{(i)} Solving $f_i = \pr_i \cmps g$ for $g$: for any $i$ in $I$ and $s$ in $S$, $f_i s = (\pr_i \cmps g)\,s = \pr_i(g\,s) = (g\,s)\,i$. By definition, $(f\trsp s) i = f_i s$, so $g = f\trsp$.
{\em Example}\/: let $f := \pr$, so $\pr_i = \pr_i \cmps \pr\trsp$.  Note that $(\pr\trsp t) i = \pr_i t = t_i = (\ide_{\Cart T} t) i$.  In fact, $\pr\trsp = \ide_{\Cart T}$.
\pbfni{(ii)} Letting $f := \gam$ in (i), $\gam_i = \pr_i \cmps \gam\trsp$.  Letting $f := \pr$ in (ii), $\pr_i = \gam_i \cmps h$.  Hence $\gam_i = \gam_i \cmps h \cmps \gam\trsp$ and $\pr_i = \pr_i \cmps \gam\trsp \cmps h$. By uniqueness, $h \cmps \gam\trsp = \ide_C$ and $\gam\trsp \cmps h = \ide_{\Cart T}$.  So $\gam\trsp$ is a bijection from $C$ to $\Cart T$ and $h$ is its inverse.
\end{therm}
This compact yet detailed proof was made possible by the {\em generic operator} $\trsp$ for {\em transposition} \cite{Bout03cgf}, defined by $(f\trsp s) i = f_i s$, as recalled inside the proof.  Dependencies may be easier to trace by writing $\mathrm{g}_f$ for $g$ in part (i) and $\mathrm{h}_{\gam,f}$ for $h$ in part (ii) of the statement.  Thus, $f := \pr$ makes $\text{h}_{\gam,\pr}$ the inverse of $\gam\trsp$.
   \par Part (i) is depicted in Figure \ref{CCPU}.
Dashed lines reflect multiple instances, one for each $i$ in $I$, forming a 3D cone.
\begin{figure}[h!]
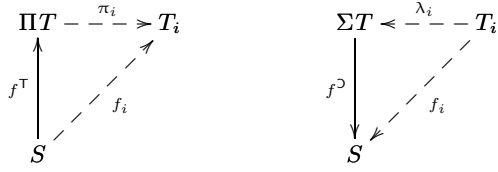

\[\bfig\small
\ptriangle/{-->}`{<-}`{<--}/[{\Cart T}`{T_i}`S;{\pr_i}`{f\trsp}`{f_i}]
\ptriangle(1200,0)/{<--}`{->}`{-->}/[{\Disu T}`{T_i}`S;{\ulam_i}`{f\unc}`{f_i}]
\efig\]
\caption{Concrete Categorical Cartesian Product and Disjoint Union}\label{CCPU}
\end{figure}
The proof of (ii) shows that every $I$-family $\gam$ with the stated property is isomorphic to $\pr$, hence ``unique up to isomorphism''.
   \par This {\em theorem} faithfully reflects the {\em definition} of products in category theory \cite[p.~19]{Pier91ctcs}, while avoiding ``unacceptable but generally accepted'' notations like $(T_i)_{i \in I}$ for just $T$ and, worse, $\langle f_i\rangle_{i \in I}$ for just $\langle f\rangle$, which was proven to be $f\trsp$.
   \par {\em Aside}\hh Clearly, $\gam \bin \Cart_{i \bin I} (C \arr T_i)$ depends on $C$ and $T$.  Similarly, $\pr$ depends on $T$, which is why some authors write $\pr^T$.  One can also see the $\pr_i$ as {\em operators} (Lamport \cite[p.~69]{Lamp03spec}) which, unlike functions, do not have a domain. A third view is offered by theories allowing for the set (or {\em class}) of all pairs, such as the algebra of relations in Tarski~\cite{Tars87stwv} or set theory with a universal set~$\mathcal{U}$ in Holmes, where axioms \cite[p.~30]{Holmxxstus} assert the existence of the {\em equality relation} $\{(x, x) \ba x \bin \mathcal{U}\}$ and {\em projection relations} such as $\{((x, y), x) \ba x \bin \mathcal{U}; y \bin \mathcal{U}\}$.
\paragraph{Sum}
This topic illustrates how objects with quite different point-wise definitions resemble each other in point-free form, more specifically as duals.
\begin{defin}[Disjoint union]\label{disu}
The {\em disjoint union} $\Disu T$ of a family $T$ of sets is the set of (ordered) pairs such that $(i, x) \membr \Disu T$ iff $i \membr \dom T$ and $x \membr T_i$.
\end{defin}
This concept is used to express choice in formal language semantics~\cite[p.~74]{Meye91itpl}.
   \par To obtain the usual point-free characterization of $\Disu T$,  we let $I := \dom T$ and define an $I$-family $\ulam$ of {\em labeling functions} by $\ulam_i \membr T_i \arr \Disu T$ and $\ulam_i x = (i, x)$.
\begin{therm}\label{sumthm} Let $T$ be an $I$-family of sets.  
\pbfni{\em(i)} For any set $S$ and $I$-family $f$ of functions with $f_i \membr T_i \arr S$, there exists a unique $g : \Disu T \arr S$ satisfying $f_i = g \cmps \ulam_i$. Specifically, $g = f\unc$ {\em (uncurry)}.
\pbfni{\em(ii)} Let $\udel$ be an $I$-family of functions with $\udel_i \membr T_i \arr D$ and the property that, for any set $S$ and $I$-family $f$ of functions with $f_i \membr T_i \arr S$, there exists a unique $h : D \arr S$ satisfying $f_i = h \cmps \udel_i$. Then there is a bijection between $D$ and $\Disu T$.
\em
\pbfni{Proof}\hh \textbf{(i)} Solving $f_i = g \cmps \ulam_i$ for $g$: for any $i$ in $I$ and $s$ in $S$, $f_i s = (g \cmps \ulam_i)\,s = g(\ulam_i s) = g(i, s)$. By definition, $f\unc (i, s) = f_i s$, so $g = f\unc$.
{\em Example}\/: let $f := \ulam$, so $\ulam_i = \ulam\unc \cmps \ulam_i$.  Note that $\ulam\unc(i, x) = \ulam_i x = (i, x) = \ide_{\Disu T} (i, x)$.   In fact, $\ulam\unc = \ide_{\Disu T}$.
\pbfni{(ii)} Letting $f := \udel$ in (i), $\udel_i = \udel\unc \cmps \ulam_i$.  Letting $f := \ulam$ in (ii), $\ulam_i = h \cmps \udel_i$.  Hence $\udel_i = \udel\unc \cmps h \cmps \udel_i$ and $\ulam_i = h \cmps \udel\unc \cmps \ulam_i$. By uniqueness, $\udel\unc \cmps h = \ide_D$ and $h \cmps \udel\unc = \ide_{\Disu T}$.  So $\udel\unc$ is a bijection from $D$ to $\Disu T$ and $h$ is its inverse.
\end{therm}
\paragraph{Relations, point-free style}
Relational properties of the Cartesian product can be similarly expressed and explored in point-free style.  For instance, given an $I$-family $R$ of relations $R_i$ from $S$ to $T_i$, let us aim for a relation $G$ from $S$ to $\Cart T$ such that $\pr_i \cmps G = R_i$.  This requires $\pr\cnv_i \cmps \pr_i \cmps G = \pr\cnv_i \cmps R_i$ and, since $\ide_{\Cart T} \sbset \pr\cnv_i \cmps \pr_i$, also $G \sbset \pr\cnv_i \cmps R_i$ (for all $i$ in $I$), hence $G \sbset \Insec i \bin I \ab \pr\cnv_i \cmps R_i$.  Let us define $\langle R\rangle$ as $\Insec i \bin I \ab \pr\cnv_i \cmps R_i$, so $t\langle R\rangle s \siff I = \mpty \sor (t \membr \Cart T \sand \evry i \bin I \ab t_i\,R_i\,s)$.
   \par In general, this only yields $\pr_i \cmps \langle R\rangle \sbset R_i$, as in category theory.  Still, if all $R_i$ share a common domain, a sharper result is $\pr_i \cmps \langle R\rangle = R_i$, the design goal.  Even sharper is the pointwise form $y\,(\pr_i \cmps \langle R\rangle)\,s \siff y\,R_i\,s \sand s \membr \Insec i \bin I \ab \dom R_i$.
   \par Assuming $I \neq \mpty$, let $T$ and $T'$ be $I$-families of sets and $R$ an $I$-family of relations $R_i$ from $T_i$ to $T'_i$. Define $\Pll R$ from $\Cart T$ to $\Cart T'$ by $\Pll R = \langle i \bin I \ab R_i \cmps \pr_i\rangle$, generalizing \cite[p.~114]{Bird97aofp}.  Now $t' (\Pll R) t \siff t' \membr \Cart T' \sand t \membr \Cart T \sand \evry i \bin I \ab t'_i\,R_i\,t_i$.  
\paragraph{Functors}
A {\em functor} in category theory maps objects to objects and arrows to arrows. Often one uses the same symbol for both maps, e.g., writing $\cart$ for $\pll$, but that is not necessary \cite[p.~30]{Bird97aofp}. Although notational economy is commendable, in ``working mathematics'' one prefers combining various concepts as convenient, yet still avoid notational confusion.  Hence, if one views relations as sets, the need to distinguish between $R \pll S$ and $R \cart S$ is evident. The operator $\pll$ is called {\em shuffle} in some functional languages and {\em parallel} by Meyer \cite[p.~36]{Meye91itpl}.
\paragraph{Tabulations}
Tabulations are helpful in point-free reasoning \cite{Bird97aofp}.  Let $R \bin X \lrar Y$. Now $f : R \arr X$ with $f(x, y) = x$ and $g : R \arr Y$ with $g(x, y) = y$ satisfy both $R = g \cmps f\cnv$ and $(f\cnv \cmps f) \insec (g\cnv \cmps g) = \ide_R$. 
   \par In category theory, a {\em tabulation} of a correspondence arrow $r$ from $X$ to $Y$ is a pair of application arrows $f$ from $Z$ to $X$ and $g$ from $Z$ to $Y$ jointly satisfying $r = g \cmps f\cnv$ and $(f\cnv \cmps f) \insec (g\cnv \cmps g) = \ide_Z$.  Similar concepts in a more general setting are discussed in Tarski \cite[p.~96]{Tars87stwv}.
\paragraph{Currying}
In its most basic form, Currying transforms a function $f \bin X \cart Y \arr Z$ into a function $f\cur$ of type $X \arr (Y \arr Z)$ such that $f\cur x\,y = f(x, y)$.  For a point-free characterization \cite[p.~72]{Bird97aofp}\cite[p.~33]{Pier91ctcs}, one uses an {\em evaluation} or {\em application} function $\alpha_{Y,Z} \bin (Y \arr Z) \cart Y \arr Z$ defined by $\alpha_{Y,Z}(g, y) = g(y)$.  One can verify that $f = \alpha_{Y,Z} \cmps (f\cur \pll \ide_Y)$; in fact, $h = f\cur \equiv \alpha_{Y,Z} \cmps (h \pll \ide_Y) = f$ (uniqueness).
\[\bfig
\morphism(0,500)|m|[X \cart Y`Z;f]
\morphism(0,500)|m|<0,-500>[{}`(Y \arr Z) \cart Y\phantom{(Y \arr )};\phantom {x}f\cur \pll \ide_Y]
\morphism(80,80)|m|<400,400>[{}`{};{\text{\raisebox{-1.5ex}{\hp3$\alpha_{Y,Z}$}}}] 
\efig\]
The categorical view is that the object $X \arr Y$ (often written $Y^X$) is fully characterized by the existence of an arrow $\alpha_{Y,Z} \bin (Y \arr Z) \cart Y \arr Z$ such that there is a unique arrow $f\cur$ satisfying $f = \alpha_{Y,Z} \cmps (f\cur \pll \ide_Y)$.
   \par In the wider setting of common mathematics, $f\cur$ is defined for any function~$f$ whose domain is a set $R$ of pairs as follows: $\tdom f\cur = \dom R$ and, for any $x$ in $\tdom f\cur$, $f\cur x \membr \{y \bin \ran R \ba (x, y) \membr R\} \arr \tran f$ with $f\cur x\,y = f(x,y)$.  {\em Example}: if $\tdom f = \{(0,5), (3, 1), (3,2)\}$ then $\tdom f\cur = \{0, 3\}$, $\dom (f\cur 3) = \{1, 2\}$ and $f\cur 3\,2 = f(3, 2)$.
   \par Conversely, uncurrying $F\unc$ is defined for any family of functions $F$ as follows: $F\unc \membr \{(x, y) \bin \tdom F \cart \Union(\tdom \cmps F) \ba y \membr \dom (F\,x)\} \arr \Union(\tran \cmps F)$ with $F\unc(x, y) = F\,x\,y$.
   \par Note that $(f\cur)\unc = f$, but $(F\unc)\cur = F$ only if $\tran F$ contains no empty function.
%
%
\section{Concluding remarks}
%
%
\noindent The first part of this paper argued that mathematical definitions and notations are the result of {\em design} and hence benefit from engineering principles. 
   \par This was illustrated for some very basic concepts that have fallen into disrepair, apparently due to habits that noted mathematicians consider ``unacceptable'' or evidence of ``glaring perversity''.  Indeed, only perversity explains the neglect for so many impeccable accounts that have been around since decades.
   \par Simple design guidelines go a long way.  In particular, definitions require diligent engineering. The discipline of always explicitly justifying design decisions incites more thoughtfulness.  Halmos~\cite{Halm60nast} provides an unsurpassed example.
   \par Judicious use of symbolism is another invaluable tool for sanity checking.  Defining concepts not just in prose but also by very simple formulas, as in Suppes \cite[p.~57, 86]{Supp72axst}, arguably would have avoided the {\em codomain} blunder.
   \par One empirical engineering principle that does {\em not} hold in mathematics is that performance requires complexity.  A ``crystal radio'' needs just one (passive) semiconductor; a software-defined radio (SDR) performs better but uses billions of transistors.  For relations and functions, the standard variant is the simplest, yet was found more versatile and general than variants with codomains.  Clearly a very onerous burden of justification rests on such variant definitions!
   \par The second part focused on making category theory practical.  Historically, practical use was not a concern, so the need for re-engineering is no surprise.
   \par Appreciating this fact took a scholar like Jos\'e Oliveira, supported by his  background in electrical engineering.  One of his recent papers \cite{Oliv14prah} points out how relation algebra requires adequate {\em preparation} to suit the intended purpose.
   \par Whereas \cite{Oliv14prah} is more advanced, we have shown how liberating category theory from unjustified constraints (similar to codomains) extends its practical applicability while preserving its essential ideas, and how full compatibility with ``ordinary mathematics'' makes the style of expression and reasoning accessible to anyone with a high school background and a taste for algebraic elegance.
   \par In this way, the basic ideas may become useful and appealing to a wider scientific community rather than remain jargon for a small group of aficionados.
\paragraph{Acknowledgement} I like to thank Jeremy Gibbons and an anonymous reviewer for many useful comments during the many revisions of this paper and their exemplary patience in discussing critical issues.
%
%
\section*{References}
%
%

\vfill\eject
\end{document}